\documentclass[12pt, reqno]{amsart}
\usepackage{amsmath,amsfonts,amssymb,amsmath,latexsym,mathrsfs}
\usepackage[all]{xy}
\usepackage{enumitem}
\usepackage{setspace}

\newcommand\Q{{\mathbb Q}}

\newcommand\Z{{\mathbb Z}}

\newtheorem{theorem}{Theorem}[section]

\newtheorem{proposition}[theorem]{Proposition}
\newtheorem{lemma}[theorem]{Lemma}

\theoremstyle{definition}
   
\newtheorem{remark}[theorem]{Remark}      
       
 \newtheorem{conjecture}[theorem]{Conjecture}

\begin{document}
\title{Homology of depth-graded motivic Lie algebras and koszulity}

\author{B. Enriquez}
\address{B.E.: IRMA (CNRS) et D\'epartement de math\'ematiques, Universit\'e de Strasbourg, 7 rue Ren\'e Descartes, 67000 
Strasbourg,  France \newline 
E-mail address: {\tt  b.enriquez@math.unistra.fr}}
\author{P. Lochak}
\address{P.L.:  CNRS and Institut Math\'ematique de Jussieu, 
Universit\'e P. et M. Curie, 4 place Jussieu, 75252 Paris Cedex 05, France\newline 
E-mail address: {\tt  pierre.lochak@imj-prg.fr}}

\begin{abstract}
The Broadhurst-Kreimer (BK) conjecture describes the Hilbert series of a bigraded Lie algebra $\mathfrak{a}$ 
related to the multizeta values. Brown proposed a conjectural description of the homology of this Lie 
algebra (homological conjecture (HC)), and showed it implies the BK conjecture. We show that a part of 
HC is equivalent to a presentation of $\mathfrak{a}$, and that the remaining part of HC
is equivalent to a weaker statement. Finally, we prove that granted the first part of HC, the remaining part of HC is equivalent to 
either of the following equivalent statements: (a) the vanishing of the third homology group of 
a Lie algebra with quadratic presentation, constructed out of the period polynomials of modular forms; (b)
the koszulity of the enveloping algebra of this Lie algebra.  
\end{abstract}

\maketitle

\section*{Introduction}

\subsection{The background}
The multizeta values (MZVs) are a family of real numbers $\zeta(n_1,..,n_s)$, where $s\geq 1$ and $n_1,\ldots,n_{s-1}\geq 1$, $n_s\geq 2$
(\cite{Z}). Denote by $\mathcal{Z}\subset \mathbb{R}$ the $\mathbb{Q}$-vector subspace spanned by these numbers; it is a subring of 
$\mathbb{R}$. Two types of algebraic relations between these numbers  are known:

$\bullet$ {\it associator relations} arise from the fact that the MZVs are the coefficients of an analytic object called the Knizhnik-Zamolodchikov
associator, and that this object satisfies the ``pentagon and duality'' relations; additional relations relating the MZVs and the complex 
number $2\pi\mathrm{i}$ can be derived in the same way from the ``hexagon'' relation (\cite{Dr,LM});

$\bullet$ {\it double shuffle and regularization relations} arise from a combinatorial study of the MZVs
 (\cite{Rac}).

The relation between the various ``associator'' relations has been elucidated in \cite{Fur1}. It can be explained 
as follows. Let $\mathcal{Z}_{assoc}$ (resp., $\tilde{\mathcal{Z}}_{assoc}$) be the $\mathbb{Q}$-ring generated by 
formal analogues $\zeta^{f}(n_1,\ldots,n_s)$ of the genuine MZVs (resp., these analogues and a formal analogue $(2\pi\mathrm{i})^f$ 
of $2\pi\mathrm{i}$), subject to the ``pentagon and duality'' relations (resp., to these relations and the ``hexagon'' relation). 
Then $\tilde{\mathcal{Z}}_{assoc}$ is isomorphic to the quadratic extension of $\mathcal{Z}_{assoc}$ generated to $(2\pi\mathrm{i})^f$ 
subject to the relation $((2\pi\mathrm{i})^f)^2=-24\zeta^f(2)$, so $\tilde{\mathcal{Z}}_{assoc}$ is a free $\mathcal{Z}_{assoc}$-module with 
basis $(1,(2\pi\mathrm{i})^f)$. 

It has also been shown (\cite{Fur2}) that the ``associator'' relations on MZVs imply the ``double shuffle'' ones.
More precisely, let $\mathcal{Z}_{ds}$ be the $\mathbb{Q}$-rings generated by the $\zeta^{f}(n_1,\ldots,n_s)$, 
subject to the double shuffle and regularization relations.
Then there is a surjective ring morphism $\mathcal{Z}_{ds}\to\mathcal{Z}_{assoc}$, taking each generator
to its analogue. 

A ring $\mathcal{Z}_{mot}$ of motivic analogues MZVs has been constructed (\cite{G1,G2}); it is linearly spanned by motivic analogues 
$\zeta^{\mathfrak{m}}(n_1,\ldots,n_s)$ of the MZVs. The motivic MZVs satisfy the associator relations (see \cite{André}, thm. 25.9.2.1), 
so there is a ring morphism $\mathcal{Z}_{assoc}\to\mathcal{Z}_{mot}$, given by $\zeta^f(n_1,\ldots,n_s)\mapsto
\zeta^{\mathfrak{m}}(n_1,\ldots,n_s)$. (By composition of the morphisms $\mathcal{Z}_{ds}\to\mathcal{Z}_{assoc}$ and 
$\mathcal{Z}_{assoc}\to\mathcal{Z}_{mot}$, one obtains a morphism $\mathcal{Z}_{ds}\to\mathcal{Z}_{mot}$, 
which had been directly constructed in \cite{Sou} before the construction of the morphism $\mathcal{Z}_{ds}\to\mathcal{Z}_{assoc}$ 
in \cite{Fur2}.) Finally, there is a (period) evaluation morphism $\mathcal{Z}_{mot}\to\mathcal{Z}$, given by $\zeta^{\mathfrak{m}}(n_1,\ldots,n_s)
\mapsto\zeta(n_1,\ldots,n_s)$. All this gives rise to a sequence of surjective ring morphisms
\begin{equation}\label{seq:Z}
\mathcal{Z}_{ds}\to\mathcal{Z}_{assoc}\to\mathcal{Z}_{mot}\to\mathcal{Z}.  
\end{equation}
The rings $\mathcal{Z}_{mot}$,  $\mathcal{Z}_{assoc}$, $\mathcal{Z}_{ds}$ are equipped with a grading, called the 
{\it weight grading,} for which $\zeta^{\mathfrak{m}}(n_1,\ldots,n_s)$ and $\zeta^{f}(n_1,\ldots,n_s)$ have weight 
$n_1+\cdots+n_s$; under the {\it direct sum conjecture,} this also defines a grading on the ring $\mathcal{Z}$. 
The rings  $\mathcal{Z}_{mot},\ldots,\mathcal{Z}$ are also equipped with an increasing filtration, the {\it depth filtration,} 
whose $d$th part is the linear span of all the elements corresponding to $(n_1,\ldots,n_s)$, with $s\leq d$. This filtration is 
compatible with the grading, unconditionally in the case of $\mathcal{Z}_{mot},\ldots,\mathcal{Z}_{ds}$, and under the direct sum 
conjecture in the case of $\mathcal{Z}$. The Hilbert series of  $\mathcal{Z}_{mot}$ with respect to the weight grading has been computed in 
\cite{B1}. It is generally conjectured that the maps from (\ref{seq:Z}) are isomorphisms compatible with the
gradings and filtrations. Under the direct sum conjecture, the depth-graded of $\mathcal{Z}$ is a bigraded algebra; 
a formula for its double Hilbert series is conjectured in \cite{BK}; we will call this the Broadhurst-Kreimer (BK) conjecture. 
Combining the BK conjecture with that of the isomorphism of $\mathcal{Z}_{ds}$,  $\mathcal{Z}_{assoc}$, $\mathcal{Z}_{mot}$
and $\mathcal{Z}$, one obtains three variants of the BK conjecture predicting the double Hilbert series of the depth-graded of 
$\mathcal{Z}_{ds}$,  $\mathcal{Z}_{assoc}$ and $\mathcal{Z}_{mot}$. 

The graded and filtered algebras $\mathcal{Z}_{ds}$,  $\mathcal{Z}_{assoc}$ and $\mathcal{Z}_{mot}$ are related to prounipotent Lie 
groups as follows. Set $\mathcal{Z}_{mot}^0:=\mathcal{Z}_{mot}/(\zeta^{\mathfrak{m}}(2))$, then there is an isomorphism of graded 
and filtered algebras $\mathcal{Z}_{mot}\simeq\mathcal{Z}_{mot}^0\otimes\mathbb{Q}[X]$, where $X$ has weight 2 all the powers 
$X,X^2,\ldots$ have depth 1 (\cite{B1}, (2.13), with the notation $\mathcal{H}\sim\mathcal{Z}_{mot}$, $\mathcal{A}\sim\mathcal{Z}_{mot}^0$; 
$X$ corresponds to $\zeta^{\mathfrak{m}}(2)$). Moreover, $\mathcal{Z}_{mot}^0$ identifies with the function algebra of a prounipotent 
group scheme $U^{\mathfrak{m}}$ arising in the theory of mixed Tate motives (\cite{B1}, Section 2.1 and main result). Its Lie algebra 
$\mathfrak{u^m}$ is graded and equipped with a decreasing filtration, which induces an isomorphism 
$\mathcal{Z}_{mot}^0\simeq U(\mathfrak{u^m})^\vee$ of graded and filtered algebras. In the same way, the fact that the rings 
$\mathcal{Z}_{assoc},\mathcal{Z}_{ds}$ are polynomial in infinitely many variables (\cite{André}, Remark 25.9.3.2) implies that there are 
isomorphisms of graded and filtered algebras $\mathcal{Z}_{assoc}\simeq\mathcal{Z}_{assoc}^0\otimes\mathbb{Q}[X]$, 
$\mathcal{Z}_{ds}\simeq\mathcal{Z}_{ds}^0\otimes\mathbb{Q}[X]$, where $X$ is a formal variable of weight 2 and such that $X,X^2,\ldots$
have depth 1, and 
$\mathcal{Z}_{assoc}^0:=\mathcal{Z}_{assoc}/(\zeta^f(2))$, $\mathcal{Z}_{ds}^0:=\mathcal{Z}_{ds}/(\zeta^f(2))$; moreover, there are 
isomorphisms of graded and filtered algebras $\mathcal{Z}_{assoc}^0\simeq U(\mathfrak{grt}_1)^\vee$ (\cite{Dr}, Prop. 5.9),  $\mathcal{Z}_{ds}^0\simeq U(\mathfrak{ds}_0)^\vee$ (\cite{Rac}, Thm. I), where $\mathfrak{grt}_1$ is the ``graded version of the Grothendieck-Teichm\"uller Lie algebra'' arising in the study of 
associators, and $\mathfrak{ds}_0$ is the ``double shuffle Lie algebra'' arising in the study of combinatorial relations between MZVs. 
All this yields compatible isomorphisms of graded and filtered algebras
$$
\mathcal{Z}_{assoc}\simeq U(\mathfrak{grt}_1)^\vee\otimes\mathbb{Q}[X], \quad 
\mathcal{Z}_{ds}\simeq U(\mathfrak{ds}_0)^\vee\otimes\mathbb{Q}[X], \quad
\mathcal{Z}_{mot}\simeq U(\mathfrak{u^m})^\vee\otimes\mathbb{Q}[X], 
$$
which in their turn yield compatible isomorphisms of bigraded algebras
$$
\mathrm{gr}_{dpth}(\mathcal{Z}_{assoc})\simeq \mathrm{gr}_{dpth}(U(\mathfrak{grt}_1))^\vee\otimes\mathbb{Q}[X], \quad 
\mathrm{gr}_{dpth}(\mathcal{Z}_{ds})\simeq \mathrm{gr}_{dpth}(U(\mathfrak{ds}_0))^\vee\otimes\mathbb{Q}[X], 
$$
$$
\mathrm{gr}_{dpth}(\mathcal{Z}_{mot})\simeq \mathrm{gr}_{dpth}(U(\mathfrak{u^m}))^\vee\otimes\mathbb{Q}[X].  
$$
The variants of the BK conjecture predicting the double Hilbert series of $\mathrm{gr}_{dpth}(\mathcal{Z}_{assoc})$ and its analogues
therefore translate into conjectural formulas for the double Hilbert series of $\mathrm{gr}_{dpth}(U(\mathfrak{grt}_1))$ and its 
analogues, which will be recalled in Section \ref{sec3} (Conjecture \ref{conj:BK}).  

Recently, F. Brown proposed a conjecture describing the homology of these bigraded Lie algebras (\cite{B2}) and implying these conjectural 
formulas. 

\subsection{Contents of the paper}
The main purpose of this paper is to partially reduce and derive some consequences of Brown's homological conjecture. 

In Section \ref{sec1}, we introduce several objects: the Lie algebras $\mathfrak{u^m},\mathfrak{grt}_1,\mathfrak{ds}_0$;  the sequence of 
morphisms $\mathfrak{u^m}\to\mathfrak{grt}_1\to\mathfrak{ds}_0$ of Lie algebras dual to the sequence of morphisms 
$\mathcal{Z}_{ds}^0\to\mathcal{Z}_{assoc}^0\to\mathcal{Z}_{mot}^0$ induced by (\ref{seq:Z}); the depth-graded version of these Lie algebras and 
of this sequence of morphisms; and an upper bound Lie algebra $\mathfrak{ls}$ of $\mathrm{gr}_{ds}(\mathfrak{ds}_0)$. In Section \ref{sec2}, we 
present some known results on these Lie algebras. In Section \ref{sec3}, we present the BK conjectures and its variants, as well as Brown's 
homological conjecture, which presents itself in four versions relative to the Lie algebras $\mathfrak{grt}_1$, $\mathfrak{ds}_0$, 
$\mathfrak{u^m}$, $\mathfrak{ls}$; we prove that the version relative to $\mathfrak{ls}$ implies all the other versions.

The main results of the paper are contained in the remaining sections. In Section \ref{sec4}, we show that a part of the homological 
conjecture for one of these Lie algebras (more precisely, the part predicting the values of the first and second homology groups) is equivalent to
a presentation of the same Lie algebra. The proof of this result is close to the proof that any positively graded Lie algebra $L$ has a presentation 
with generating space $H_1(L)$ and relation space $H_2(L)$, where $H_i(L)$ denotes the $i$th Lie algebra homology group of $L$ with trivial 
coefficients (\cite{H}, Section 3); the arguments from \cite{H} are themselves analogues of 
those of \cite{S}, Chap. 2 in the pro-$p$ group situation. In Section \ref{sec5}, we show that the remaining part of the  homological conjecture, 
more precisely, the vanishing of the homology groups of order $\geq 3$, is equivalent 
to a weaker statement, namely the vanishing of the third homology group. Again, the technical result used here is an analogue of a homological 
result for pro-$p$ groups from \cite{S}. In Section \ref{sec6}, we relate the homological conjecture with the structure of a Lie algebra
$\mathfrak{M}_0$ with quadratic presentation, constructed out of the period polynomials of cusp modular forms. We show that 
granted the first part of the homological conjecture relative to any of the Lie algebras $\mathfrak{grt}_1,\ldots,\mathfrak{ls}$, 
the remaining part of this conjecture is equivalent to either of the following equivalent statements: (a) the vanishing of the Lie 
algebra homology group $H_3(\mathfrak{M}_0,\mathbf{k})$; (b) the koszulity of the enveloping algebra $U(\mathfrak{M}_0)$. 

\section{Motivic and related Lie algebras} \label{sec1}

In this section, we recall the construction of a Lie algebra $\mathfrak{g^m}$ attached to the 
category of mixed Tate motives (Subsection \ref{subsect:tate}); the realization of its Lie subalgebra $\mathfrak{u^m}$ by derivations of 
a free Lie algebra in two generators (Subsections \ref{subsect:12} and \ref{subsect:13}); and the construction of Lie algebras 
related to $\mathfrak{u^m}$ (Subsection \ref{subsect:14}). 

\subsection{Motivic background} \label{subsect:tate}

The Tannakian category $\mathrm{MT}(\mathbb{Z})$ of mixed Tate motives over $\mathbb{Z}$ can be defined unconditionally 
(see \cite{DG}; for a survey 
see \cite{André}). It is equipped with the de Rham fiber functor. Let $\mathcal{G}_{\mathrm{MT}(\mathbb{Z})}$
be the automorphism group of this functor; this is a group scheme over $\mathbb{Q}$. It decomposes as a semidirect product
$\mathcal{G}_{\mathrm{MT}(\mathbb{Z})}=\mathcal{U}_{\mathrm{MT}(\mathbb{Z})}\rtimes\mathbb{G}_m$, 
where $\mathcal{U}_{\mathrm{MT}(\mathbb{Z})}$ is a prounipotent $\mathbb{Q}$-group scheme. Its Lie algebra
decomposes as $\mathfrak{g^m}=\mathfrak{u^m}\rtimes\mathbb{Q}$. It follows from work of Borel and Beilinson 
that there exist elements $s_{2i+1}$, $i\geq 1$ of $\mathfrak{u^m}$, freely generating $\mathfrak{u^m}$, 
and a right inverse $\mathbb{G}_m\to\mathcal{G}_{\mathrm{MT}(\mathbb{Z})}$ of the projection 
$\mathcal{G}_{\mathrm{MT}(\mathbb{Z})}\to\mathbb{G}_m$, for which each $s_{2i+1}$ has weight $2i+1$. 
We will henceforth view $\mathfrak{u^m}$ as graded by this grading. 

\subsection{The Poisson-Ihara Lie algebra $(\mathfrak{g},\langle,\rangle)$} \label{subsect:12}


Let $e_0,e_1$ be free noncommutative variables of weight $1$, and let $\mathfrak{g}:=\mathbb{L}(e_0,e_1)$ be the free 
$\mathbb{Q}$-Lie algebra generated by these variables; its bracket will be denoted $[,]$. 
There is a linear map $\mathfrak{g}\to\mathrm{Der}(\mathfrak{g})$, given by $f\mapsto D_f$, where $D_f:e_0\mapsto0$, 
$e_1\mapsto[e_1,f]$. When equipped with the {\it Poisson-Ihara Lie bracket $\langle,\rangle$} given by 
$$
\langle f,g\rangle:=[f,g]+D_f(g)-D_g(f), 
$$ 
$\mathfrak{g}$ is a Lie algebra, and the map $f\mapsto D_f$ is a Lie algebra morphism. Moreover, $(\mathfrak{g},\langle,\rangle)$
is graded for the weight grading. The Lie algebra $(\mathfrak{g},[,])$ is also equipped with a grading, for which $e_1$ has degree 1 and $e_0$
has degree 0. This equips $(\mathfrak{g},\langle,\rangle)$ with a grading, the {\it depth grading}. The {\it depth filtration} on 
$(\mathfrak{g},\langle,\rangle)$ is the decreasing Lie algebra filtration given by $F^i(\mathfrak{g}):=\oplus_{j|j\geq i}
\{\mathrm{part\ of\ }\mathfrak{g}\mathrm{\ of\ depth\ }j\}$. 

\subsection{The Lie algebra morphism $\mathfrak{u^m}\to(\mathfrak{g},\langle,\rangle)$} \label{subsect:13}

The motivic Galois group $\mathcal{G}_{\mathrm{MT}(\mathbb{Z})}$ acts naturally on the fundamental groupoids of certain 
geometric objects of $\mathrm{MT}(\mathbb{Z})$. In the case of $\mathbb{P}^1_{\mathbb{Q}}-\{0,1,\infty\}$, this induces a 
weight-graded Lie algebra morphism $\mathfrak{u^m}\to\mathrm{Der}(\mathfrak{g})$. It can be shown that this morphism
factors through the inclusion $(\mathfrak{g},\langle,\rangle)\subset\mathrm{Der}(\mathfrak{g})$ and therefore gives rise to a
graded Lie algebra morphism $\mathfrak{u^m}\to(\mathfrak{g},\langle,\rangle)$. Moreover, it follows from \cite{B1} that this morphism 
is {\it injective.} 

\subsection{Lie algebras related to $\mathfrak{u^m}$} \label{subsect:14} 

In \cite{Dr,Rac}, explicit weight-graded Lie subalgebras $\mathfrak{grt}_1$ and $\mathfrak{ds}_0$ 
of $(\mathfrak{g},\langle,\rangle)$ are introduced. 

There is a sequence of Lie algebra morphisms
\begin{equation}\label{seq:LAs}
\mathfrak{u^m}\simeq\mathrm{im}(\mathfrak{u^m}\to\mathfrak{g})\hookrightarrow\mathfrak{grt}_1
\hookrightarrow\mathfrak{ds}_0\subset(\mathfrak{g},\langle,\rangle)
\end{equation}
where as mentioned above, the initial isomorphism follows from \cite{B1}, and the next injections follow from work of Drinfeld and 
Ihara, and from \cite{Fur2}, respectively. 

\subsection{The Lie algebra $\mathfrak{ls}$}

The Lie algebra $\mathfrak{ds}_0$ is defined to be the set of elements of the free algebra generated by $e_0,e_1$ which are primitive 
for two coproducts, the shuffle and the stuffle 
coproduct. The associated graded of the stuffle coproduct for the depth filtration can be computed explicitly. This gives rise to an 
explicit upper bound Lie algebra for the depth-graded of the double shuffle Lie algebra $\mathfrak{ds}_0$, called the ``linearized shuffle''
Lie algebra and denoted $\mathfrak{ls}$. We have therefore a double inclusion $\mathrm{gr}_{dpth}(\mathfrak{ds}_0)\subset
\mathfrak{ls}
\subset(\mathfrak{g},\langle,\rangle)$ (see \cite{B2}). 

\begin{remark}
This space was introduced earlier in \cite{IKZ}, and in dual form, in \cite{G1,G2}; moreover, the Lie algebras $\mathfrak{ls}$ and 
$\mathfrak{ds}_0$ are respectively isomorphic to the Lie algebras introduced in \cite{E1,E2} under the notation 
$\mathsf{ARI}_{\mathrm{ent}}^{\mathrm{al/al}}$ and $\mathsf{ARI}_{\mathrm{ent}}^{\mathrm{al/il}}$.
\end{remark}


\subsection{Bigraded Lie algebras}

The depth filtration of the Lie algebra $(\mathfrak{g},\langle,\rangle)$ induces a filtration on each of its weight-graded Lie subalgebras; 
the associated graded Lie algebra is then a (depth,weight)-bigraded Lie subalgebra of $\mathrm{gr}_{dpth}(\mathfrak{g},
\langle,\rangle)=(\mathfrak{g},\langle,\rangle)$. The sequence (\ref{seq:LAs}) therefore gives rise to a sequence 
\begin{equation} \label{seq:gr:LAs}
\mathrm{gr}_{dpth}(\mathfrak{u^m})\hookrightarrow\mathrm{gr}_{dpth}(\mathfrak{grt}_1)
\hookrightarrow\mathrm{gr}_{dpth}(\mathfrak{ds}_0)\hookrightarrow\mathfrak{ls}\subset (\mathfrak{g},\langle,\rangle)
\end{equation}
Numerical experimentation indicates that the first two injections of 
(\ref{seq:LAs}), and the first three injections of (\ref{seq:gr:LAs}) could in fact be isomorphisms. 

\section{Known results on depth-graded Lie algebras}\label{sec2}

In this section, we recall some results on the Lie algebras $\mathfrak{a}=\mathrm{gr}_{dpth}(\mathfrak{u^m}),\ldots,\mathfrak{ls}$ from 
(\ref{seq:gr:LAs}). We first introduce the space $\mathsf{P}$ of even period polynomials (Subsection \ref{subsect:31}). Using this space, 
we make explicit a system of generators and relations for the Lie algebras $\mathfrak{a}$ of (\ref{seq:gr:LAs}) (Subsection 
\ref{subsect:32}). We draw some consequences for $\mathfrak{ls}$ (Subsection \ref{subsect:33}) and for the other Lie algebras 
$\mathfrak{a}$ (Subsection \ref{subsect:34}). 

\subsection{The spaces $\mathsf{S}$ and $\mathsf{P}$} \label{subsect:31}

Let $\mathsf{S}$ denote the (complex) vector space of cusp forms for the full modular group $\mathrm{PSL}_2(\mathbb{Z})$,
which decomposes as $\mathsf{S} =\bigoplus_n \mathsf{S}_{2n} $ where $\mathsf{S}_{2n}$ denotes the space of forms
of weight $2n$ ($n\ge 0$). One sets $\mathbb{S}(s):=\sum_{n\geq 0}\mathrm{dim}(\mathsf{S}_{2n})s^{2n}$. Then  
\begin{equation}\label{def:S}
 \mathbb{S}(s)={s^{12} \over (1-s^4)(1-s^6)}.
\end{equation}
Let $\tilde{\mathsf{P}}_{2n}$ denote the $\mathbb{Q}$-vector space of all even period polynomials of degree $2n-2$, 
which are the polynomials in $\Q[X,Y]$ satisfying certain symmetry conditions: evenness in $X$ and $Y$, antisymmetry 
with respect to the exchange of $X$ and $Y$ and a functional equation. The additional condition of divisibility 
by $X^2Y^2$ defines a hyperplane $\mathsf{P}_{2n}\subset\tilde{\mathsf{P}}_{2n}$. The Eichler-Shimura correspondence sets up 
a linear isomorphism $\mathsf{S}_{2n}\simeq\mathsf{P}_{2n}\otimes_{\mathbb{Q}}\mathbb{C}$ for any $n\geq 0$. We set 
$\mathsf{P}:=\bigoplus_n \mathsf{P}_{2n}$. 

\subsection{Generators and relations in $\mathfrak{ls}$}\label{subsect:32}

\subsubsection{Generators}\label{subsect:gens}

The depth one subspace $\mathfrak{ls}_1$ of $\mathfrak{ls}$ is 1-dimensional in each odd weight $\geq 3$ and 0-dimensional 
in all the other weights, therefore if $\mathbb{O}(s):=\sum_{n\geq 0} \mathrm{dim}(\mathfrak{ls}_1[n])s^n$, then 
\begin{equation}\label{def:O}
\mathbb{O}(s)={s^3 \over 1-s^2}
\end{equation}
(where $[n]$ means the part of weight $n$). On the other hand, there exists an injective linear map 
$\mathsf{e}:\mathsf{P}\to\mathfrak{ls}_4$, compatible with the weight gradings on both sides (see \cite{B2}).  

\subsubsection{Relations}\label{subsect:map}

As mentioned above, there is an isomorphism of graded vector spaces $\mathfrak{ls}_1\simeq
X^3\mathbb{C}[X^2]$, where the degree on the right-hand side is the degree in $X$. On the other hand, there is an 
injection $\mathsf{P}\hookrightarrow X^2Y^2\mathbb{C}[X^2,Y^2]^{as}$, where $as$ means antisymmetry in $X^2,Y^2$, 
and where $\mathsf{P}_{2n}$ maps to the part of total degree $2n-2$. The composition 
$$
\mathsf{P}\hookrightarrow X^2Y^2\mathbb{C}[X^2,Y^2]^{as}\stackrel{XY\times-}{\to}
X^3Y^3\mathbb{C}[X^2,Y^2]^{as}\simeq\Lambda^2(X^3\mathbb{C}[X^2])\simeq\Lambda^2(\mathfrak{ls}_1)
$$
is then an injective graded linear map $\mathsf{P}\to\Lambda^2(\mathfrak{ls}_1)$. On the other hand, the Lie bracket 
of $\mathfrak{ls}$ induces a linear map $\Lambda^2(\mathfrak{ls}_1)\to\mathfrak{ls}_2$. It follows from \cite{GKZ,Sch}
that these two maps combine into an exact sequence 
$$
0\to\mathsf{P}\to\Lambda^2(\mathfrak{ls}_1)\to\mathfrak{ls}_2\to 0. 
$$

\subsection{Generators and relations in the Lie algebras $\mathrm{gr}_{dpth}(\mathfrak{u^m}),\ldots,\mathfrak{ls}$}
\label{subsect:gens:0}\label{subsect:map:0}

Let $\mathfrak{a}$ be one of the Lie algebras $\mathrm{gr}_{dpth}(\mathfrak{u^m}),\ldots,\mathfrak{ls}$ in sequence (\ref{seq:gr:LAs}). 
There are isomorphisms $\mathfrak{a}_i\simeq\mathfrak{ls}_i$ for $i=1,2$ between the depth $1,2$ subspaces of $\mathfrak{a}$
and $\mathfrak{ls}$. Therefore there is an exact sequence 
$$
0\to\mathsf{P}\to\Lambda^2(\mathfrak{a}_1)\to\mathfrak{a}_2\to 0. 
$$
isomorphic to the analogous exact sequence for the Lie algebra $\mathfrak{ls}$.  

\subsection{Consequences for $\mathfrak{ls}$}\label{subsect:33}

\subsubsection{A model Lie algebra $\mathfrak{M}$} We define a (weight, depth)-bigraded Lie algebra $\mathfrak{M}$, 
which will be called the ``model'' Lie algebra, by the following presentation. The space of generators is the direct sum of: 

$\bullet$ a copy of the (weight, depth)-bigraded vector space $\mathfrak{ls}_1$; 

$\bullet$ a copy $\mathsf{P}\{4\}$ of the space $\mathsf{P}$, where $\{4\}$ means that for each $n$, the weight $n$ component of 
$\mathsf{P}$ is placed in bidegree $(n,4)$. 

The map $\mathsf{P}\to\Lambda^2(\mathfrak{ls}_1)$ constructed in Subsection \ref{subsect:map} induces a bigraded linear map 
$\mathsf{P}\{2\}\hookrightarrow\Lambda^2(\mathfrak{ls}_1)$, where $\{2\}$ has the same meaning as above. We then obtain a 
composed map 
$$
\mathsf{P}\{2\}\hookrightarrow\Lambda^2(\mathfrak{ls}_1)\simeq\mathbb{L}_2(\mathfrak{ls}_1)
\hookrightarrow\mathbb{L}_2(\mathfrak{ls}_1\oplus\mathsf{P}\{4\})\hookrightarrow\mathbb{L}(\mathfrak{ls}_1\oplus\mathsf{P}\{4\}), 
$$
where $\mathbb{L}(V)$ is the free Lie algebra generated by a (possibly graded) vector space $V$, and $\mathbb{L}_i(V)$
is its degree $i$ part (with respect to $V$). The space of defining relations of $\mathfrak{M}$ is defined to be the image of this map. 
We have therefore 
$$
\mathfrak{M}:=\mathbb{L}(\mathfrak{ls}_1\oplus\mathsf{P}\{4\})/(\mathsf{P}\{2\}), 
$$
where $(\mathsf{P}\{2\})$ is the ideal generated by the image of the map $\mathsf{P}\{2\}\hookrightarrow
\mathbb{L}(\mathfrak{ls}_1\oplus\mathsf{P}\{4\})$. 

\subsubsection{A morphism $\mathfrak{M}\to\mathfrak{ls}$}

It follows from Subsections \ref{subsect:gens} and \ref{subsect:map} that there exists a morphism of bigraded Lie algebras
\begin{equation}\label{morph:M:ls}
\varphi:\mathfrak{M}\to\mathfrak{ls},
\end{equation}
defined by restriction to the space of generators of $\mathfrak{M}$ as follows: the restriction of $\varphi$ to $\mathfrak{ls}_1$
is the canonical injection $\mathfrak{ls}_1\hookrightarrow\mathfrak{ls}$; its restriction to $\mathsf{P}\{4\}$
is the composed map $\mathsf{P}\{4\}\simeq\mathsf{P}\stackrel{\mathsf{e}}{\to}\mathfrak{ls}_4\hookrightarrow\mathfrak{ls}$. 

\subsubsection{The homology of $\mathfrak{ls}$}

If $\mathfrak{a}$ is a Lie algebra, we denote by $H_\cdot(\mathfrak{a})$ its homology with coefficients in the trivial module $\mathbf{k}$,
given by the homology of the complex $\ldots\to\Lambda^3(\mathfrak{a})\stackrel{[,]\otimes\mathrm{id}}{\to}
\Lambda^2(\mathfrak{a})\stackrel{[,]}{\to}\mathfrak{a}\to 0$. This is a graded cocommutative coalgebra, depending 
functorially on $\mathfrak{a}$. 

It follows from the definition of $\mathfrak{M}$ that 
$$
H_1(\mathfrak{M})\simeq\mathfrak{ls}_1\oplus\mathsf{P}\{4\}. 
$$
Moreover, the composed map $\mathsf{P}\{2\}\to\Lambda^2(\mathfrak{ls}_1)\to\Lambda^2(\mathfrak{M})
\stackrel{[,]}{\to}\mathfrak{M}$ is zero. This induces a map 
$$
\mathsf{P}\{2\}\to H_2(\mathfrak{M}). 
$$
One then checks that these maps are compatible with coproducts, so that the diagram
\begin{equation}\label{Pls}
\xymatrix{\mathsf{P}\{2\} \ar[r]\ar[d]& \Lambda^2(\mathfrak{ls}_1)\ar[r]& \Lambda^2(\mathfrak{ls}_1\oplus\mathsf{P}\{4\})\ar[d]
\\ H_2(\mathfrak{M})\ar[rr]& & \Lambda^2(H_1(\mathfrak{M}))}
\end{equation}
commutes. Using the morphism $H_\cdot(\mathfrak{M})\to H_\cdot(\mathfrak{ls})$ induced by the 
morphism (\ref{morph:M:ls}), one obtains: 

\begin{proposition}
\begin{enumerate}
\item[i)] There exists a bigraded linear map $\mathfrak{ls}_1\oplus\mathsf{P}\{4\}\to H_1(\mathfrak{ls})$; 
\item[ii)] there exists a bigraded linear map $\mathsf{P}\{2\}\to H_2(\mathfrak{ls})$; 
\item[iii)] the diagram (\ref{Pls}) commutes. 
\end{enumerate}
\end{proposition}

\subsection{Consequences for the Lie algebras $\mathrm{gr}_{dpth}(\mathfrak{u^m}),\ldots,\mathfrak{ls}$}\label{subsect:34}

Let $\mathfrak{a}$ be one of the Lie algebras $\mathrm{gr}_{dpth}(\mathfrak{u^m}),\ldots,\mathfrak{ls}$ in sequence (\ref{seq:gr:LAs}). 

\subsubsection{A Lie algebra $\mathfrak{M}_0$} We set 
$$
\mathfrak{M}_0:=\mathbb{L}(\mathfrak{ls}_1)/(\mathsf{P}\{2\}). 
$$
Then $\mathfrak{M}_0$ is a (weight, depth)-bigraded Lie algebra. Then $\mathfrak{M}$ can be identified with the 
free product of Lie algebras $\mathfrak{M}_0*\mathbb{L}(\mathsf{P}\{4\})$. 

\subsubsection{A morphism $\mathfrak{M}_0\to\mathfrak{a}$}

It follows from Subsection \ref{subsect:map:0} that there exists a morphism of bigraded Lie algebras
$$
\varphi_0:\mathfrak{M}_0\to\mathfrak{a},
$$
defined by restriction to the space of generators of $\mathfrak{M}_0$ as follows: the restriction of $\varphi_0$ to $\mathfrak{ls}_1$
is the canonical injection $\mathfrak{ls}_1\simeq\mathfrak{a}_1\hookrightarrow\mathfrak{a}$. 

\subsubsection{The homology of $\mathfrak{a}$}

The map $\mathfrak{ls}_1\to\mathfrak{M}_0$ induces a linear map $\mathfrak{ls}_1\to H_1(\mathfrak{M}_0)$. 
The linear map $\mathsf{P}\{2\}\to \Lambda^2(\mathfrak{ls}_1)$ is such that the composite map 
$\mathsf{P}\{2\}\to\Lambda^2(\mathfrak{ls}_1)\hookrightarrow\Lambda^2(\mathfrak{M}_0)\stackrel{[,]}{\to}\mathfrak{M}_0$
is zero.  One derives from there a linear map $\mathsf{P}\{2\}\to H_2(\mathfrak{M}_0)$. One checks that the diagram 
\begin{equation}\label{diag:PM0}
\xymatrix{ \mathsf{P}\{2\}\ar[r]\ar[d]&\Lambda^2(\mathfrak{ls}_1)\ar[d] \\ H_2(\mathfrak{M}_0)\ar[r]& \Lambda^2(H_1(\mathfrak{M}_0))}
\end{equation}
commutes. Using functoriality of Lie algebra homology with respect to the morphism $\varphi_0:\mathfrak{M}_0\to\mathfrak{a}$, and 
the fact that the morphism $\mathfrak{M}_0\to\mathfrak{a}$ is an isomorphism in depths 1 and 2, one obtains: 

\begin{proposition}\label{prop:comm:diag}\label{prop:32}
\begin{enumerate}
\item[i)] There exists a bigraded linear map $\mathfrak{ls}_1\to H_1(\mathfrak{a})$, which is an isomorphism in depth $1$; 
\item[ii)] there exists a bigraded linear map $\mathsf{P}\{2\}\to H_2(\mathfrak{a})$, which is an isomorphism in depth $2$; 
\item[iii)] the diagram (\ref{diag:PM0}) commutes. 
\end{enumerate}
\end{proposition}

\section{Conjectures on depth-graded Lie algebras}\label{sec3}

We present the motivic/combinatorial analogue of the Broadhurst-Kreimer conjecture, which is a statement about Hilbert series (Subsection 
\ref{subsect:hilb}). We then present the homological counterpart of this conjecture, due to Brown, in Subsection 
\ref{subsect:homology}. We then discuss the relations between both conjectures (Subsection \ref{subsection:relation}). 

\subsection{Hilbert series} \label{subsect:hilb}

Based on numerical experimentation, the authors of \cite{BK} formulated a conjecture on dimensions of the depth-graded spaces of 
(real) multizeta values. The authors of \cite{B2,CGS} then proposed motivic/combinatorial analogues of this conjecture, which we now 
recall. 

\begin{conjecture} \label{conj:BK} (\cite{BK,B2,CGS})
{\it 
Let $\mathfrak{a}$ be any of the Lie algebras $\mathrm{gr}_{dpth}(\mathfrak{u^m}),\ldots,\mathfrak{ls}$ 
in sequence (\ref{seq:gr:LAs}). Let $S(\mathfrak{a})$ be its symmetric algebra and let $S(\mathfrak{a})[n,d]$ denote its piece of weight $n$ and
depth $d$. Then there is an identity of formal series in two variables $s, t$:
$$
\sum_{n,d\geq 0} \mathrm{dim}S(\mathfrak{a})[n,d]\cdot s^nt^d=
{1\over  1-\mathbb{O}(s)t+\mathbb{S}(s)t^2-\mathbb{S}(s)t^4}, 
$$
where $\mathbb{O}(s)$, $\mathbb{S}(s)$ are defined in (\ref{def:O}), (\ref{def:S}).} 
\end{conjecture}
Conjecture \ref{conj:BK} for the Lie algebra $\mathfrak{a}$ will be denoted $\mathrm{HS}(\mathfrak{a})$
(``Hilbert series statement for $\mathfrak{a}$''). Then: 

\begin{proposition} There holds the sequence of implications 
$$\mathrm{HS}(\mathfrak{ls})\Rightarrow
\mathrm{HS}(\mathrm{gr}_{dpth}(\mathfrak{ds}_0))\Rightarrow\mathrm{HS}(\mathrm{gr}_{dpth}(\mathfrak{grt}_1))\Rightarrow
\mathrm{HS}(\mathrm{gr}_{dpth}(\mathfrak{g^m})).
$$ 
\end{proposition}

{\em Proof.} Set $\mathfrak{a}_0:=\mathfrak{g^m}$, $\mathfrak{a}_1:=\mathrm{gr}_{dpth}(\mathfrak{grt}_1)$,
$\mathfrak{a}_2:=\mathrm{gr}_{dpth}(\mathfrak{ds}_0))$, $\mathfrak{a}_3:=\mathfrak{ls}$. For $i=0,\ldots,3$, let 
$f_i(s,t):=\sum_{n,d\geq 0}\mathrm{dim}S(\mathfrak{a}_i)[n,d]\cdot s^nt^d\in\mathbb{Z}_{\geq 0}[t][[s]]$. The inclusions
$\mathfrak{a}_0\subset\cdots\subset\mathfrak{a}_3$ imply $f_0(s,t)\leq\cdots\leq f_3(s,t)$, where $\leq$
means that the difference belongs to $\mathbb{Z}_{\geq 0}[t][[s]]$. The structure result on $\mathfrak{g^m}$ (see Subsection \ref{subsect:tate}) 
implies that $f_0(s,1)=1/(1-\mathbb{O}(s))$. 

Fix $i\in\{0,1,2\}$. The statement  $\mathrm{HS}(\mathfrak{a}_{i+1})$ implies that $f_{i+1}(s,t)=1/(1-\mathbb{O}(s)t+\mathbb{S}
(s)t^2-\mathbb{S}(s)t^4)=:g(s,t)$. It follows that 
\begin{equation}\label{ineq:1}
f_i(s,t)\leq g(s,t).
\end{equation} 
On the other hand, 
\begin{equation}\label{ineq:2}
f_i(s,1)\geq f_0(s,1)=g(s,1).
\end{equation}  
If $f,g$ belong to $\mathbb{Z}_{\geq 0}[t]$ and are such that $f(t)\leq g(t)$ and $f(1)\geq g(1)$, then $f=g$. 
It follows that if $F,G$ belong to $\mathbb{Z}_{\geq 0}[t][[s]]$ and are such that $F(s,t)\leq G(s,t)$ and $F(s,1)\geq G(s,1)$, 
then $F=G$. (\ref{ineq:1}) and (\ref{ineq:2}) then imply $f_i(s,t)=g(s,t)$, therefore $\mathrm{HS}(\mathfrak{a}_i)$. 
\hfill \qed\medskip 

\subsection{Homology}\label{subsect:homology} 

\begin{conjecture} \label{conj:HC} (\cite{B2}) {\it Let $\mathfrak{a}$ be any of the Lie algebras $\mathrm{gr}_{dpth}(\mathfrak{u^m}),\ldots,\mathfrak{ls}$ in sequence (\ref{seq:gr:LAs}).
\begin{enumerate}
\item[i)] There exists a bigraded linear isomorphism $\mathfrak{ls}_1\oplus\mathsf{P}\{4\}\simeq H_1(\mathfrak{a})$, such that the diagram 
$$
\xymatrix{\mathfrak{ls}_1\oplus\mathsf{P}\{4\}\ar[r]\ar[rd] & H_1(\mathfrak{a})\ar[d]\\ &H_1(\mathfrak{ls})}
$$
commutes; 
\item[ii)] the map $\mathsf{P}\{2\}\to H_2(\mathfrak{a})$ is a linear isomorphism; 
\item[iii)] for any $k\geq 3$, $H_k(\mathfrak{a})=0$. 
\end{enumerate}
}\end{conjecture}
Conjecture \ref{conj:HC} for the Lie algebra $\mathfrak{a}$ will be denoted $\mathrm{HC}(\mathfrak{a})$ (``homology conjecture 
for Lie algebra $\mathfrak{a}$''). Taking into account Proposition \ref{prop:32}, ii), HC$(\mathfrak{a})$, ii), says that $H_2(\mathfrak{a})$ is concentrated in depth 2.  
It also follows from Proposition \ref{prop:comm:diag} that for any Lie algebra $\mathfrak{a}$ in sequence (\ref{seq:gr:LAs}), the diagram 
$$
\xymatrix{\mathsf{P}\{2\} \ar[r]\ar[d]& \Lambda^2(\mathfrak{ls}_1)\ar[d]
\\ H_2(\mathfrak{a})\ar[r]& \Lambda^2(H_1(\mathfrak{a}))}
$$ 
commutes. 

\subsection{Relation between the Hilbert series and the homological conjectures} \label{subsection:relation}

It has been proved in \cite{B2} that for $\mathfrak{a}$ one of the Lie algebras $\mathrm{gr}_{dpth}(\mathfrak{u^m}),\ldots,\mathfrak{ls}$ 
in sequence (\ref{seq:gr:LAs}), HC($\mathfrak{a}$) implies HS($\mathfrak{a}$). 

\section{Homological conjecture and presentation}\label{sec4}

In this section, we establish the equivalence between a part of the homological conjecture \ref{conj:HC} and presentation results 
for the corresponding Lie algebras. 
In Subsection \ref{subsect:51}, we show one implication, and in Subsection \ref{subsect:52}, we show the converse implication. 
As mentioned in the Introduction, the results of this section are close to those of \cite{H} and \cite{S}, Chap. 2. 

\subsection{From presentation to homology} \label{subsect:51}

If $V$ is a vector space and $I\subset\mathbb{L}(V)$ is a Lie ideal, then the 
Lie analogue of Hopf's formula (see \cite{H}, Proposition 5.6) gives
$$H_2(\mathbb{L}(V)/I)=I/[\mathbb{L}(V),I];$$
it can be derived from the definition of the homology of $\mathbb{L}(V)/I$ using a long exact sequence. 

If now $R$ is a vector subspace of $\mathbb{L}_2(V)$ and $I:=(R)\subset\mathbb{L}(V)$ is the ideal generated by $R$, 
then $I$ decomposes as $I=R\oplus[V,R]\oplus[V,[V,R]]\oplus\cdots$ according to the degree with respect to $V$, 
while $[\mathbb{L}(V),I]=[V,R]\oplus[V,[V,R]]\oplus\cdots$; it follows that 
\begin{equation}\label{hopf}
H_2(\mathbb{L}(V)/(R))=R. 
\end{equation}
One derives from there 
\begin{equation}
H_2(\mathfrak{M})=\mathsf{P}\{2\}. \label{comp:H2}
\end{equation}
One also computes
\begin{equation} \label{comp:H1}
H_1(\mathfrak{M})=\mathfrak{M}/[\mathfrak{M},\mathfrak{M}]=\mathfrak{ls}_1\oplus\mathsf{P}\{4\}. 
\end{equation}
(\ref{comp:H1}) and (\ref{comp:H2}) then imply: 

\begin{lemma} \label{lemma:1}
Let $\mathfrak{a}$ be one of the Lie algebras $\mathrm{gr}_{dpth}(\mathfrak{u^m}),\ldots,\mathfrak{ls}$ 
in sequence (\ref{seq:gr:LAs}). If there is a bigraded Lie algebra isomorphism 
$\mathfrak{M}\simeq\mathfrak{a}$, such that the diagram $\xymatrix{\mathfrak{M}\ar[r]\ar[rd] & \mathfrak{a}\ar[d]\\ &\mathfrak{ls}}$
commutes, then $\mathrm{HC}(\mathfrak{a})$ i) and ii) hold.
\end{lemma}  

\subsection{From homology to presentation} \label{subsect:52}

Let $\mathfrak{g}$ be a positively graded Lie algebra, so $\mathfrak{g}=\mathfrak{g}_1\oplus\mathfrak{g}_1\oplus\ldots$, 
with $[\mathfrak{g}_i,\mathfrak{g}_j]\subset\mathfrak{g}_{i+j}$. The grading of $\mathfrak{g}$ then induces a positive grading on 
each of the homology groups of $\mathfrak{g}$.
There is a quotient map $\mathfrak{g}\to H_1(\mathfrak{g})=\mathfrak{g}/[\mathfrak{g},\mathfrak{g}]$, which is surjective and graded. 
Set $H:=H_1(\mathfrak{g})$ and let $H=H_1\oplus H_2\oplus\ldots$ be the degree decomposition of $H$. Let us choose a graded 
right inverse $H\to\mathfrak{g}$ of the projection map $\mathfrak{g}\to H$. It gives rise to a surjective, graded Lie algebra morphism 
$\mathbb{L}(H)\to\mathfrak{g}$. Let $I$ be the kernel of this morphism. Then $I$ is a graded ideal of $\mathbb{L}(H)$; let 
$I=I_2\oplus I_3\oplus\ldots$ be the degree decomposition of $I$. Then there is a graded Lie algebra isomorphism 
$\mathbb{L}(H)/I\simeq\mathfrak{g}$. Equation (\ref{hopf}) then implies that $H_2(\mathfrak{g})=I/[\mathbb{L}(H),I]$.
 {\it Assume that $H_2(\mathfrak{g})$ is concentrated in degree 2.} Then 
$$
I_3=[H_1,I_2], \quad I_4=[H_2+[H_1,H_1],I_2]+[H_1,I_3], 
$$
$$
I_5=[H_3+[H_1,H_2]+[H_1,[H_1,H_1]],I_2]+[H_2+[H_1,H_1],I_3]+[H_1,I_4], \quad\mathrm{etc.}
$$
It follows that for any $n>0$, 
\begin{equation} \label{incl}
I_{n+2}\subset \sum_{k\geq 1}\sum_{\stackrel{n_1,\ldots,n_k>0,}{n_1+\cdots+n_k=n}}
[H_{n_1},[H_{n_2},\ldots,[H_{n_k},I_2]]]. 
\end{equation}
Since conversely, the right-hand side of (\ref{incl}) is contained in its left-hand side, one has 
for any $n>0$, $I_{n+2}=\sum_{k\geq 1}\sum_{\stackrel{n_1,\ldots,n_k>0,}{n_1+\cdots+n_k=n}}
[H_{n_1},[H_{n_2},\ldots,[H_{n_k},I_2]]]$. It follows that $I=(I_2)$. We also have $H_2(\mathfrak{g})\simeq I_2$,
and the injective composed map $H_2(\mathfrak{g})\simeq I_2\subset\mathbb{L}_2(H_1)\subset\mathbb{L}_2(H)\simeq\Lambda^2(H_1(\mathfrak{g}))$ 
can be identified with the coproduct map $H_2(\mathfrak{g})\to\Lambda^2(H_1(\mathfrak{g}))$. It follows that this coproduct map 
is injective. 

All this implies: 

\begin{lemma} \label{lemma:identif}
If $\mathfrak{g}$ is a positively graded Lie algebra, such that $H_2(\mathfrak{g})$ is concentrated in degree 2, then
the coproduct map $H_2(\mathfrak{g})\to\Lambda^2(H_1(\mathfrak{g}))$ is injective, and there exists an isomorphism of 
graded Lie algebras $\mathfrak{g}\simeq \mathbb{L}(H_1(\mathfrak{g}))/(H_2(\mathfrak{g}))$, where $H_2(\mathfrak{g})$
is viewed as a subspace of $\mathbb{L}_2(H_1(\mathfrak{g}))$ via the sequence of maps $H_2(\mathfrak{g})\to
\Lambda^2(H_1(\mathfrak{g}))\simeq\mathbb{L}_2(H_1(\mathfrak{g})))$. 
\end{lemma}

This lemma can be adapted to the case when $\mathfrak{g}$ is also equipped with a second grading, compatible with the first one. 
Let $\mathfrak{a}$ be any of the Lie algebras $\mathrm{gr}_{dpth}(\mathfrak{u^m}),\ldots,\mathfrak{ls}$ in sequence (\ref{seq:gr:LAs})
and assume that it satisfies $\mathrm{HC}(\mathfrak{a})$, i) and ii). Then applying the graded version of Lemma \ref{lemma:identif} 
to $\mathfrak{a}$ equipped with the depth grading as the first grading, and the weight grading as the second grading, 
we obtain a bigraded isomorphism $\mathfrak{M}\simeq\mathfrak{a}$. 

\subsection{\ } Combining Lemma \ref{lemma:identif} with Lemma \ref{lemma:1}, we obtain: 

\begin{theorem} \label{thm53}
For $\mathfrak{a}$ one of the Lie algebras $\mathrm{gr}_{dpth}(\mathfrak{u^m}),\ldots,\mathfrak{ls}$ 
in sequence (\ref{seq:gr:LAs}), the conjunction of $\mathrm{HC}(\mathfrak{a})$ i) and ii) is equivalent to the existence of a bigraded Lie algebra isomorphism $\mathfrak{M}\simeq
\mathfrak{a}$, such that the diagram $$
\xymatrix{\mathfrak{M}\ar[r]\ar[rd] & \mathfrak{a}\ar[d]\\ &\mathfrak{ls}}
$$
commutes. 
\end{theorem}

\begin{remark} Under $\mathrm{HC}(\mathfrak{a})$ i) and ii), Theorem \ref{thm53} gives a presentation of $\mathfrak{a}$. 
When $\mathfrak{a}=\mathrm{gr}_{dpth}(\mathfrak{ds}_0)$, this presentation was first conjectured in \cite{E1}, \S 17, 18
and \cite{E2}, \S 7,8. 
\end{remark}

\section{Reduction of the homological conjecture}\label{sec5}
 
Let $\mathfrak{g}=\oplus_{i>0}\mathfrak{g}_i$ be a positively graded Lie algebra over a field $\mathbf{k}$ of characteristic 0.  
We assume that the graded pieces $\mathfrak{g}_i$ are finite dimensional. 

Let $\mathcal{M}od$ be the category whose objects are the graded $\mathfrak{g}$-modules $V$ of the form 
$V=\bigoplus_{i\in\mathbb{Z}}  V_i$, where each graded piece $V_i$ is finite dimensional and $V_i=0$ for $i\ll 0$, and the morphisms 
are the graded $\mathfrak{g}$-module morphisms. Let $\mathcal{V}ec$ denote the category of graded vector spaces with 
the same conditions. 

Define the {\it valuation} $v$ as the map from the set of objects of any of these categories to 
$\Z\cup\{\infty\}$ given by $v(V)= {\rm min} \{ i\in \Z, V_i\ne 0 \}$ for $V$ a nonzero object and $v(0)=\infty$. 
For $a$ an integer, define the {\it shift operator by $a$}, denoted $V\mapsto V[a]$, as the self-map of the set of objects 
of any of the categories $\mathcal{M}od$ or $\mathcal{V}ec$, given by $V[a]_i=V_{i+a}$ for any $i\in\mathbb{Z}$. 

For $k\ge 0$, there is a functor $H_k: {\mathcal Mod}\rightarrow {\mathcal Vec}$ given by $H_k(V):=H_k(\mathfrak{g},V)$
(Lie algebra homology group); we call it the {\it $k$-th homology functor.} 

Let $n>0$ be an integer. Denote again by $\mathbf{k}$ the trivial $\mathfrak{g}$-module (one-dimensional, concentrated in degree 0), 
and {\it assume that $H_n(\mathfrak{g},\mathbf{k})=0$.} 

Let $V$ be a nonzero object in $\mathcal{M}od$; set $v:=v(V)$. For $s\geq 0$, define the {\it truncated module} 
$V\{s\}$ as the object of $\mathcal{M}od$ given by $V\{s\}_i=V_i$ for $i\geq v+s$, $V\{s\}_i=0$ otherwise (in particular, 
$V\{0\}=V$). Then for any $s\geq 0$, there is an exact sequence in $\mathcal{M}od$
$$
0\to V\{s+1\}\to V\{s\}\to\mathbf{k}[s+v]^{\mathrm{dim}(V_{s+v})}\to0, 
$$
 where the last module is a direct sum of copies of shifts of the trivial module $\mathbf{k}$. The assumption $H_n(\mathbf{k})=0$
implies $H_n(\mathbf{k}[s+v])=0$, therefore the homology long exact sequence implies that for any $s\geq 0$, the map 
$H_n(V\{s+1\})\to H_n( V\{s\})$ is {\it onto.} For any $s\geq 0$, each of the maps of the sequence $H_n( V\{s\})\to H_n( V\{s-1\})\to
\cdots\to H_n(V\{0\})=H_n(V)$ is therefore onto, so that
\begin{equation}\label{star}
\mathrm{\ the\ composed\ map\ }H_n( V\{s\})\to H_n(V)\mathrm{\ is\ onto.}
\end{equation}
On the other hand, $H_n(V\{s\})$ is the homology group of a complex constructed out of the graded vector space 
$\Lambda^\cdot(\mathfrak{g})\otimes V\{s\}$ , which is nonzero only in degrees $\geq v+s$, therefore 
\begin{equation}\label{starstar}
H_n(V\{s\})\mathrm{\ is\ nonzero\ only\ in\ degrees\ }\geq v+s. 
\end{equation}
For any integer $k\in\mathbb{Z}$ and any $s\geq \mathrm{max}(0,k-v+1)$, (\ref{starstar}) implies that (degree $k$ part of 
$H_n(V\{s\})$)=0 and (\ref{star}) implies that the map (degree $k$ part of $H_n(V\{s\}))\to$(degree $k$ part of 
$H_n(V)$) is onto, so that (degree $k$ part of $H_n(V))=0$. It follows that $H_n(V)=0$.  Therefore 
\begin{equation}\label{Hn:zero}
\mathrm{for\ any\ object\ }V{\ of\ }\mathcal{M}od, H_n(V)=0. 
\end{equation}
Let $U(\mathfrak{g})$ be the universal enveloping algebra of $\mathfrak{g}$. When equipped with the grading induced by $\mathfrak{g}$
and with the action of $\mathfrak{g}$ by left multiplication, $U(\mathfrak{g})$ is an object of $\mathcal Mod$. The truncation 
$U(\mathfrak{g})\{1\}$ can be identified with the kernel $U(\mathfrak{g})_+$ of the counit map $U(\mathfrak{g})\to\mathbf{k}$. Then 
there is an exact sequence in $\mathcal Mod$
$$
0\to U(\mathfrak{g})_+\to U(\mathfrak{g})\to\mathbf{k}\to 0. 
$$ 
Shapiro`s lemma implies that $H_i(\mathfrak{g},U(\mathfrak{g}))=0$ for any integer $i>0$. When $i=0$, 
$H_0(\mathfrak{g},U(\mathfrak{g}))=\simeq\mathbf{k}\simeq H_0(\mathfrak{g},\mathbf{k})$. The homology long exact sequence 
then implies an isomorphism $H_{i+1}(\mathfrak{g},\mathbf{k})\simeq H_i(\mathfrak{g},U(\mathfrak{g})_+)$ for any integer $i$. 
Combining this equality for $i=n$ with the specialization of (\ref{Hn:zero}) for $V=U(\mathfrak{g})_+$, we obtain 
$H_{n+1}(\mathbf{k})=0$. 

We summarize these results as follows: 

\begin{proposition}\label{prop:nilp}
Let $\mathfrak{g}=\oplus_{i>0}\mathfrak{g}_i$ be a positively graded Lie algebra over a field $\mathbf{k}$ of characteristic $0$, 
such that all the graded pieces $\mathfrak{g}_i$ are finite dimensional. Let $n$ be an integer $\geq 0$. If $H_n(\mathfrak{g},\mathbf{k})=0$, 
then for any $m\geq n$, $H_m(\mathfrak{g},\mathbf{k})=0$. 
\end{proposition}

\begin{remark} Proposition \ref{prop:nilp} illustrates a general principle in homological algebra for (pro)nilpotent objects;
compare in particular Proposition 21 in \cite{S} (\S 4) for pro-$p$ groups. \hfill\qed
\end{remark}

An immediate corollary of Proposition \ref{prop:nilp} is the following result.  

\begin{theorem} \label{h3}
For $\mathfrak{a}$ one of the Lie algebras $\mathrm{gr}_{dpth}(\mathfrak{u^m}),\ldots,\mathfrak{ls}$ in sequence
(\ref{seq:gr:LAs}), item iii) of the Homological Conjecture \ref{conj:HC} is equivalent to $H_3(\mathfrak{a})=0$.  
\end{theorem}

\section{Homological conjecture and koszulity}\label{sec6}
 
In this section, we show that for any Lie algebra $\mathfrak{a}=\mathrm{gr}_{dpth}(\mathfrak{ds}_0),\ldots,\mathfrak{ls}$ in sequence
(\ref{seq:gr:LAs}), the conjecture $\mathrm{HC}(\mathfrak{a})$ implies two equivalent statements, namely the vanishing of the homology group 
$H_3(\mathfrak{M}_0,\mathbf{k})$, and the koszulity of the algebra $U(\mathfrak{M}_0)$, equipped with the 
depth grading. This section is organized as follows. In Subsection \ref{subsect:71}, we compute the homology of free products
of algebras. In Subsection \ref{subsect:72}, we use this result to show that $\mathrm{HC}(\mathfrak{a})$ implies the vanishing of 
$H_3(\mathfrak{M}_0,\mathbf{k})$. In Subsection \ref{subsect:73}, we use a result of Goncharov (\cite{G1}) to prove an equality 
$(\mathfrak{ls}_1\otimes\mathsf{P})\cap(\mathsf{P}\otimes\mathfrak{ls}_1)=0$ on the generators and relations of $\mathfrak{M}_0$.
After recalling some results on Koszul algebras (Subsection \ref{subsect:koszul}), we study Lie algebras with a quadratic 
presentation in Subsection \ref{subsect:74}: we prove that for such a Lie algebra with space of
generators $V$ and space of relations $R\subset\mathbb{L}_2(V)$, such that $(V\otimes R)\cap(R\otimes V)=0$, the vanishing of the third Lie algebra homology group is equivalent to the koszulity of its enveloping algebra. In Subsection \ref{subsect:75}, we gather all these results
to prove the main theorem \ref{main:thm:koszul}.   

\subsection{Homology of free products of algebras}\label{subsect:71}

Let $\mathbf{k}$ be a field. The category of $\mathbf{k}$-algebras with unit is equipped with a coproduct operation
(see \cite{Bbk}, Chap. III, p. 195, exercice 6). For $A,B$ two $\mathbf{k}$-algebras with unit, we denote by $A*B$ their 
coproduct in this category (also called the {\it free product} of $A$ and $B$). 

The coproduct property implies that if $A,B$ admit presentations of the form $A=T(V)/(R)$, $B=T(W)/(S)$, 
where $V,W$ are vector spaces and $R,S$ are vector subspaces of $T(V),T(W)$, then a presentation of $A*B$ is 
\begin{equation}\label{copdt:presented}
A*B\simeq T(V\oplus W)/(R\oplus S), 
\end{equation}
where we use the canonical injections $R\subset T(V)\hookrightarrow T(V\oplus W)$, $S\subset T(W)\hookrightarrow T(V\oplus W)$. 

It is proved in {\it loc. cit.} that $A*B$ is equipped with a complete increasing filtration 
$P_0\subset P_1\subset\cdots\subset A*B$, where $P_n$ is the image of the part of tensor degree $\leq n$
of the tensor algebra $T(A\oplus B)$ under the algebra morphism $T(A\oplus B)\to A*B$ induced by the linear map 
$A\oplus B\to A*B$, direct sum of the canonical maps $A\to A*B$, $B\to A*B$, and that there exist isomorphisms 
$$
P_0\simeq \mathbf{k}, \quad P_ {2n}/P_{2n-1}\simeq ((A/\mathbf{k})\otimes(B/\mathbf{k}))^{\otimes n}\oplus ((B/\mathbf{k})\otimes(A/\mathbf{k}))^{\otimes n} \mathrm{\ if\ }n>0, 
$$
$$ 
P_{2n+1}/P_{2n}\simeq  ((A/\mathbf{k})\otimes(B/\mathbf{k}))^{\otimes n}\otimes (A/\mathbf{k})\oplus ((B/\mathbf{k})\otimes(A/\mathbf{k}))^{\otimes n}\otimes(B/\mathbf{k}) \mathrm{\ if\ }n\geq 0. 
$$
One derives from there that if $A,B$ are equipped with augmentation morphisms $A\to\mathbf{k}$, $B\to\mathbf{k}$ with kernels $I,J$, then 
injection followed by product induces an isomorphism 
\begin{equation}\label{map:copdt}
\mathbf{k}\oplus\big(\bigoplus_{n>0}(I\otimes J)^{\otimes n}\oplus(J\otimes I)^{\otimes n}\big)
\oplus\big(\bigoplus_{n\geq 0}(I\otimes J)^{\otimes n}\otimes I\oplus(J\otimes I)^{\otimes n}\otimes J\big)
\stackrel{\sim}{\to}A*B. 
\end{equation}
Moreover, the augmentation morphisms $A\to\mathbf{k}$, $B\to\mathbf{k}$ induce an augmentation morphism 
$A*B\to\mathbf{k}$, whose kernel $\mathcal I$ can be identified with the image of the sum of all the summands of (\ref{map:copdt})
except $\mathbf{k}$. Together with the fact that (\ref{map:copdt}) is an isomorphism, this implies: 

\begin{lemma}\label{lemma:71}
The map $(A*B)\otimes_A I\oplus (A*B)\otimes_B J\to{\mathcal I}$ induced by the product is an isomorphism of left 
$A*B$-modules. 
\end{lemma}

Let $C\to\mathbf{k}$ be an augmented algebra and let $N$ be a $C$-module. Recall that $\mathbf{k}$ is a $C$-module through the 
augmentation map. Let $K\subset C$ be the kernel of this map; there is an exact sequence of $C$-modules 
$$
0\to K\to C\to\mathbf{k}\to 0. 
$$
As $C$ is a free $C$-module, one has $\mathrm{Tor}^C_i(C,\mathbf{k})=0$ for any $i\geq 1$. The long exact sequence for Tor then 
then yields the isomorphism 
\begin{equation}\label{iso:augm:alg}
\mathrm{Tor}_i^C(\mathbf{k},N)=\mathrm{Tor}_{i-1}^C(K,N) \quad\mathrm{for\ any\ }i\geq 1. 
\end{equation}
Let now $A\to\mathbf{k}$, $B\to\mathbf{k}$ be augmented algebras and let $M$ be an $A*B$-module. The isomorphism 
(\ref{iso:augm:alg}) implies that for any $i>1$, there is an isomorphism 
\begin{equation}\label{iso:first}
\mathrm{Tor}_i^{A*B}(\mathbf{k},M)\simeq\mathrm{Tor}_{i-1}^{A*B}(\mathcal{I},M). 
\end{equation}
Lemma \ref{lemma:71} then implies that 
\begin{equation}\label{iso:second}
\mathrm{Tor}_i^{A*B}(\mathcal{I},M)\simeq\mathrm{Tor}_{i-1}^{A*B}((A*B)\otimes_A I,M)\oplus
\mathrm{Tor}_{i-1}^{A*B}((A*B)\otimes_B J,M)
\end{equation}
According to (\ref{map:copdt}), inclusion followed by product induces an isomorphism of $A$-modules
$$
A\otimes X\stackrel{\sim}{\to}A*B, \quad\mathrm{where}\quad
X=\mathbf{k}\oplus\big(\bigoplus_{n>0}(J\otimes I)^{\otimes n}\big)\oplus\big(\bigoplus_{n>0}(J\otimes I)^{\otimes n}\otimes J\big)
$$ 
where $A\otimes X$ is viewed as a free $A$-module and $A*B$ is viewed as an $A$-module under left multiplication. It follows that 
$A*B$ is free as an $A$-module, and therefore is flat over $A$. Likewise, $A*B$ is flat over $B$. One derives from there isomorphisms
$\mathrm{Tor}_{i-1}^{A*B}((A*B)\otimes_A I,M)\simeq\mathrm{Tor}_{i-1}^A(I,M)$ and likewise exchanging $A$ and $B$, therefore
\begin{equation}\label{iso:third}
\mathrm{Tor}_{i-1}^{A*B}((A*B)\otimes_A I,M)\oplus\mathrm{Tor}_{i-1}^{A*B}((A*B)\otimes_B J,M)
\simeq \mathrm{Tor}_{i-1}^A(I,M)\oplus\mathrm{Tor}_{i-1}^B(J,M).
\end{equation}
Using again (\ref{iso:augm:alg}), we obtain 
\begin{equation}\label{iso:fourth}
\mathrm{Tor}_{i-1}^A(I,M)\oplus\mathrm{Tor}_{i-1}^B(J,M)\simeq \mathrm{Tor}_i^A(\mathbf{k},M)\oplus\mathrm{Tor}_i^B(\mathbf{k},M).
\end{equation}
Combining (\ref{iso:first}), (\ref{iso:second}), (\ref{iso:third}) and (\ref{iso:fourth}), we obtain: 
\begin{proposition}\label{prop:amalg}
If $A,B$ are augmented algebras and $M$ is any $A*B$-module, then 
$$
\mathrm{Tor}_i^{A*B}(\mathbf{k},M)=\mathrm{Tor}_i^{A}(\mathbf{k},M)\oplus \mathrm{Tor}_i^{B}(\mathbf{k},M) 
$$
for any $i\geq 2$. 
\end{proposition}


\subsection{Homological conjecture and vanishing of $H_3(\mathfrak{M}_0,\mathbf{k})$}\label{subsect:72}

Let $\mathfrak{a}$ be one of the Lie algebras $\mathrm{gr}_{dpth}(\mathfrak{ds}_0),\ldots,\mathfrak{ls}$ in sequence (\ref{seq:gr:LAs}). 
Assume that $\mathrm{HC}(\mathfrak{a})$ i) and ii) hold. By Theorem \ref{thm53}, we then have a Lie algebra isomorphism $\mathfrak{a}
\simeq\mathfrak{M}$. 
According to Theorem \ref{h3}, $\mathrm{HC}(\mathfrak{a})$, iii) says that $H_3(\mathfrak{a},\mathbf{k})=0$. Given the 
isomorphism $\mathfrak{a}\simeq\mathfrak{M}$, this last statement is equivalent to  
\begin{equation}\label{eq:0}
H_3(\mathfrak{M},\mathbf{k})=0. 
\end{equation} 
Set $A:=U(\mathfrak{M}_0)$ and $B:=T(\mathsf{P}\{4\})$. Identity (\ref{copdt:presented}) implies that $A*B$ identifies with 
$U(\mathfrak{M})$ as an augmented algebra. Proposition \ref{prop:amalg} and the identity $H_i(\mathfrak{g},M)
=\mathrm{Tor}_i^{U(\mathfrak{g})}(\mathbf{k},M)$ relating the homology of a Lie algebra $\mathfrak{g}$ with values in a 
$\mathfrak{g}$-module $M$ with the Tor groups of its enveloping algebra (\cite{W}, Cor. 7.3.6) imply the equality 
$H_i(\mathfrak{M},\mathbf{k})=H_i(\mathfrak{M}_0,\mathbf{k})\oplus H_i(\mathbb{L}(\mathsf{P}\{4\}),\mathbf{k})$
for any $i\geq 2$. The equality $H_i(\mathbb{L}(\mathsf{P}\{4\}),\mathbf{k})=0$ for $i\geq 2$ then implies that 
$H_i(\mathfrak{M},\mathbf{k})=H_i(\mathfrak{M}_0,\mathbf{k})$ for $i\geq 2$, and therefore
$H_3(\mathfrak{M},\mathbf{k})=H_3(\mathfrak{M}_0,\mathbf{k})$. Equality (\ref{eq:0}) is therefore equivalent to 
$H_3(\mathfrak{M}_0,\mathbf{k})=0$. 

We have proved: 

\begin{proposition}\label{prop:van:H3}
Assume that for $\mathfrak{a}$ one of the Lie algebras $\mathrm{gr}_{dpth}(\mathfrak{ds}_0),\ldots,\mathfrak{ls}$ in sequence 
(\ref{seq:gr:LAs}), $\mathrm{HC}(\mathfrak{a})$ i) and ii) hold. Then $\mathrm{HC}(\mathfrak{a})$ iii) is equivalent to 
$H_3(\mathfrak{M}_0,\mathbf{k})=0$. 
\end{proposition}

\subsection{The equality $(\mathfrak{ls}_1\otimes\mathsf{P})\cap(\mathsf{P}\otimes\mathfrak{ls}_1)=0$}\label{subsect:73}

The three last lines of \cite{G1} contain the computation of the homology of a complex denoted there $\mathcal{D}_{\bullet,3}\to\mathcal{D}_{\bullet,2}\otimes\mathcal{D}_{\bullet,1}\to\Lambda^3(\mathcal{D}_{\bullet,1})$: 
this complex is acyclic except at $\mathcal{D}_{\bullet,3}$. The dual complex is 
$\Lambda^3(\mathfrak{ls}_1)\stackrel{[,]\otimes\mathrm{id}}{\to}\mathfrak{ls}_2\otimes\mathfrak{ls}_1\stackrel{[,]}{\to}\mathfrak{ls}_3$; 
the result of \cite{G1} then implies that the map 
\begin{equation}\label{map:gonchy}
\Lambda^3(\mathfrak{ls}_1)\stackrel{[,]\otimes\mathrm{id}}{\to}\mathfrak{ls}_2\otimes\mathfrak{ls}_1
\end{equation} is injective. 

Since the kernel of the map $[,]:\Lambda^2(\mathfrak{ls}_1)\to\mathfrak{ls}_2$ is equal to $\mathsf{P}$, 
the kernel of the map (\ref{map:gonchy}) is equal to $(\mathsf{P}\otimes\mathfrak{ls}_1)\cap\Lambda^3(\mathfrak{ls}_1)$. 
It then follows form the injectivity of (\ref{map:gonchy}) that 
\begin{equation}\label{inters:zero}
(\mathsf{P}\otimes\mathfrak{ls}_1)\cap\Lambda^3(\mathfrak{ls}_1)=0.
\end{equation} 

The subspace $(\mathfrak{ls}_1\otimes\mathsf{P})\cap(\mathsf{P}\otimes\mathfrak{ls}_1)$ of $\mathfrak{ls}_1^{\otimes 3}$
is contained in $\mathfrak{ls}_1\otimes\Lambda^2(\mathfrak{ls}_1)$ (as $\mathsf{P}$ consists of antisymmetric tensors)
and in $\Lambda^2(\mathfrak{ls}_1)\otimes\mathfrak{ls}_1$ (for the same reason), therefore it is contained in the 
intersection of these spaces, namely $\Lambda^3(\mathfrak{ls}_1)$. It follows that $(\mathfrak{ls}_1\otimes\mathsf{P})\cap(\mathsf{P}\otimes\mathfrak{ls}_1)$ is contained in the 
intersection of $\Lambda^3(\mathfrak{ls}_1)$ with $\mathsf{P}\otimes\mathfrak{ls}_1$, which by (\ref{inters:zero})
is zero. We have proved:

\begin{proposition} \label{prop:VPPV}
The equality $(\mathfrak{ls}_1\otimes\mathsf{P})\cap(\mathsf{P}\otimes\mathfrak{ls}_1)=0$ holds 
(equality of subspaces of $\mathfrak{ls}_1^{\otimes 3}$). 
\end{proposition} 

\subsection{The notion of koszulity}\label{subsect:koszul}

Let $A=\oplus_{i\geq 0}A_i$ be a graded connected algebra over a field $\mathbf{k}$, with finite dimensional graded pieces. 
Set $A_+:=\oplus_{i>0}A_i$; equip the space $T_+(A_+):=\oplus_{s>0}A_+^{\otimes s}$ with the bidegree (syzygy, weight) by declaring 
$A_{n_1}\otimes\cdots\otimes A_{n_s}$ of bidegree $(s,n_1+\cdots+n_s)$. The space $T_+(A_+)$ is equipped with a differential 
of bidegree $(-1,0)$ given by $a_1\otimes\cdots\otimes a_s\mapsto\sum_{i=1}^{s-1}(-1)^{i+1}a_1\otimes\cdots\otimes
a_ia_{i+1}\otimes\cdots\otimes a_s$. The corresponding homology group $H(T_+(A_+),d)$ inherits a bigrading, so it 
decomposes as $H(T_+(A_+),d)=\oplus_{i,j>0}\mathrm{Tor}^A_{ij}(\mathbf{k},\mathbf{k})$. Observe that the group 
$\mathrm{Tor}^A_{ij}(\mathbf{k},\mathbf{k})$ vanishes unless $i\leq j$. 

The graded algebra $A$ is called {\it Koszul} iff the groups $\mathrm{Tor}^A_{ij}(\mathbf{k},\mathbf{k})$ vanish unless $i\neq j$
(\cite{PP}, Def. 1, p. 19, and {\it loc. cit.}, identity relating Tor and Ext, end of p. 3). 

When $A$ is quadratically presented, i.e., $A=T(V)/(R)$, where $V$ is a finite dimensional vector space and $R\subset V^{\otimes 2}$
is a vector subspace, and $V$ has degree 1, the diagonal homology $\oplus_{i>0}\mathrm{Tor}^A_{ii}(\mathbf{k},\mathbf{k})$
identifies as a graded vector space with the positive part of the graded coalgebra $C(A):=\oplus_{n\geq 0}C_n(A)$, where 
$C_n(A)=\cap_{i=0}^{n-2}V^{\otimes i}\otimes R\otimes V^{\otimes n-2-i}$ (\cite{LV}, Proposition 3.3.2). 

\subsection{Koszulity of enveloping algebras}\label{subsect:74} 

Let $V$ be a finite dimensional vector space and $R\subset\mathbb{L}_2(V)$ be a vector subspace. Let $\mathfrak{g}:=
\mathbb{L}(V)/(R)$ be the corresponding Lie algebra. Its enveloping algebra will be denoted by $A$, so 
$A=U(\mathfrak{g})=T(V)/(R)$. If we assign to $V$ the degree 1, this is a quadratic algebra. 

There are canonical isomorphisms $H_i(\mathfrak{g})\simeq C_i(A)$ for $i=1,2$. 
Together with the equality $H_i(\mathfrak{g})=\oplus_{j}\mathrm{Tor}_A^{ij}(\mathbf{k},\mathbf{k})$ (see \cite{W}, 
Cor. 7.3.6) and the isomorphism $C_i(A)\simeq \mathrm{Tor}_A^{ii}(\mathbf{k},\mathbf{k})$ (\cite{LV}, Proposition 3.3.2), 
this implies
\begin{equation}\label{partial:result}
\mathrm{Tor}_A^{ij}(\mathbf{k},\mathbf{k})=0 \mathrm{\ for\ }i\neq j \mathrm{\ and\ }i=1,2. 
\end{equation}

Assume now that $(V\otimes R)\cap(R\otimes V)=0$. Then $C_n(A)=0$ for any $n\geq 3$, therefore (using again
\cite{LV}, Proposition 3.3.2) 
\begin{equation}\label{partial:result:2}
\mathrm{Tor}_A^{ii}(\mathbf{k},\mathbf{k})=0\mathrm{\ for\ }i\geq 3. 
\end{equation}
Given (\ref{partial:result}), the koszulity of $A$ is equivalent to 
$$
\mathrm{Tor}_A^{ij}(\mathbf{k},\mathbf{k})=0 \mathrm{\ for\ }i\neq j \mathrm{\ and\ }i\geq 3. 
$$
Given (\ref{partial:result:2}), the latter statement is equivalent to 
$$
\mathrm{Tor}_A^{ij}(\mathbf{k},\mathbf{k})=0 \mathrm{\ for\ }i\geq 3\mathrm{\ and\ any\ }j>0. 
$$
This statement is equivalent to 
$$
\oplus_{j>0}\mathrm{Tor}_A^{ij}(\mathbf{k},\mathbf{k})=0\mathrm{\ for\ }i\geq 3, 
$$
which by \cite{W}, Cor. 7.3.6 is equivalent to 
$$
H_i(\mathfrak{g},\mathbf{k})=0\mathrm{\ for\ }i\geq 3.  
$$
Taking into account Proposition \ref{prop:nilp}, this is equivalent to the vanishing of $H_3(\mathfrak{g},\mathbf{k})$. 

We have proved: 

\begin{proposition} \label{prop:kosz}
Let $V$ be a finite dimensional vector space and $R\subset\mathbb{L}_2(V)$ be a vector subspace. Assume that $(V\otimes R)\cap
(R\otimes V)=0$. Let $\mathfrak{g}:=\mathbb{L}(V)/(R)$ be the Lie algebra with space of generators $V$ and space of relations $R$. 
The following conditions are equivalent: 
\begin{itemize}
\item[i)] the algebra $U(\mathfrak{g})$, equipped with the grading for which $V$ has degree $1$, is Koszul; 
\item[ii)] $H_3(\mathfrak{g},\mathbf{k})=0$. 
\end{itemize}
\end{proposition}
Note that in this proposition, the category of finite dimensional vector spaces can be replaced by that of $\mathbb{Z}_{\geq 0}$-graded 
vector spaces with finite dimensional graded pieces.  

Combining Proposition \ref{prop:VPPV} with this variant of Proposition \ref{prop:kosz} applied to $V=\mathfrak{ls}_1$, $R=\mathsf{P}$ and $\mathfrak{g}=\mathfrak{M}_0$, we obtain: 

\begin{proposition}\label{prop:76}
The following statements are equivalent: 
\begin{itemize}
\item[i)] $H_3(\mathfrak{M}_0,\mathbf{k})=0$; 
\item[ii)] the algebra $U(\mathfrak{M}_0)$, equipped with the depth grading, is Koszul. 
\end{itemize}
\end{proposition}

\begin{remark} As $\mathfrak{M}_0$ is graded by depth, so is its homology group $H_3(\mathfrak{M}_0,\mathbf{k})$. 
If $H_3(\mathfrak{M}_0,\mathbf{k})_d$ denotes the depth $d$ part of this group, then $H_3(\mathfrak{M}_0,\mathbf{k})_d=0$
for $d=1,2$ for obvious reasons. Proposition \ref{prop:VPPV} implies that $H_3(\mathfrak{M}_0,\mathbf{k})_3=0$. 
\end{remark}

\subsection{\ }\label{subsect:75}

Combining Theorem \ref{thm53}, Proposition \ref{prop:van:H3} and Proposition \ref{prop:76}, we obtain: 

\begin{theorem}\label{main:thm:koszul}
Let $\mathfrak{a}$ be one of the Lie algebras $\mathrm{gr}_{dpth}(\mathfrak{ds}_0),\ldots,\mathfrak{ls}$ in sequence 
(\ref{seq:gr:LAs}). 
\begin{enumerate}
\item[(a)] The conjunction of $\mathrm{HC}(\mathfrak{a})$ i) and ii) is equivalent to the isomorphism $\mathfrak{a}\simeq\mathfrak{M}$. 
\item[(b)] Granted the conjunction of $\mathrm{HC}(\mathfrak{a})$ i) and ii), $\mathrm{HC}(\mathfrak{a})$ iii) is equivalent to either
of the following statements: (1) the vanishing of $H_3(\mathfrak{M}_0,\mathbf{k})$; (2) the koszulity of the algebra $U(\mathfrak{M}_0)$, 
equipped with the depth grading (for which $\mathfrak{ls}_1$ has depth $1$). 
\end{enumerate}
\end{theorem}

\subsection*{Acknowledgements} This work was started during a visit of both authors at the workshop 
``Grothendieck-Teichm\"uller Theory and Multiple Zeta Values'', organized in Spring 2013 at the 
Isaac Newton Institute for Mathematical Sciences. The authors thank H. Gangl for his collaboration at an early stage of this work.

\end{document}